\documentclass[12pt,a4paper]{article}
\usepackage[utf8]{inputenc}
\usepackage{amsmath,amssymb,mathtools,amsthm}
\usepackage{todonotes}
\usepackage[margin=1in]{geometry}
\usepackage{comment}
\usepackage{graphicx}
\usepackage{hyperref}
\usepackage{cite}
\usepackage{authblk}
\usepackage{xspace}

\title{Spectral similarity for Barab\'asi-Albert and Chung-Lu models}
\author[1,2]{Adam Glos\thanks{aglos@iitis.pl}}
\affil[1]{Institute of Theoretical and Applied Informatics,\protect\\
	Polish Academy of Sciences,\protect\\
	Ba{\l}tycka 5, 44-100 Gliwice, Poland}
\affil[2]{Institute of Informatics, Silesian University of Technology,\protect\\
	ul. Akademicka 16, 44-100 Gliwice, Poland}
\date{}

\newcommand{\R}{\mathbb R}

\newcommand{\w}{\mathbf w}
\newcommand{\dd}{\mathrm d}
\newcommand{\ee}{\mathrm e}
\newcommand{\ii}{\mathrm i}

\newcommand{\CL}{Chung-Lu\xspace}
\newcommand{\BA}{Barab\'asi-Albert\xspace}
\newcommand{\ER}{Erd\H{o}s-R\'enyi\xspace}

\newcommand{\ket}[1]{\ensuremath{|#1\rangle}}
\newcommand{\bra}[1]{\ensuremath{\langle#1|}}

\newcommand{\proj}[1]{\ensuremath{\ket{#1}\bra{#1}}}

\begin{document}
\maketitle

\begin{abstract}
	In the paper we have analyzed spectral similarity between Barab\'asi-Albert and
Chung-lu models. We have shown the similarity of spectral distribution for
sufficiently large Barab\'asi-Albert parameter value. Contrary, extreme
eigenvalues and principal eigenvector are not similar for those model. We
provide applications of obtained results related to the spectral graph theory
and efficiency of quantum spatial search
\end{abstract}

\section{Introduction} 

From the very first paper concerning random graphs \cite{erdos1959random}, many
new random graph models have appeared in the literature
\cite{albert2002statistical,watts1998collective,
	chung2002connected,penrose2003random}. They can be divided into two large
classes: those for which the edges are added independently, like \ER
\cite{erdos1959random} or \CL \cite{chung2002connected}, and those for which
such independence does not appear, like \BA \cite{albert2002statistical} or
Watts-Strogatz~\cite{watts1998collective}. The first class is usually easier to
analyze, especially in the context of its spectral properties
\cite{chung2011spectra,erdHos2013spectral}. Unfortunately graph models from the
second class have usually more desirable properties, like small-world or
power-law degree distribution \cite{albert2002statistical}.

Clearly, the models where the edges are added dependently cannot be described by
any model with independently added edges. However, it may be possible that some
properties match. As an example, let us consider the degree distribution between
\ER and Watts-Strogatz models. While they are different, the degree
distributions in both models are homogeneous. However, the spectral properties
of these models in general differ. This includes spectral distribution and
second greatest eigenvalue distribution of adjacency matrix
\cite{farkas2001spectra}. This is important as, spectral properties have
application in both computer science and physics
\cite{anand2011shannon,glos2018vertices,chakraborty2016spatial,faccin2013degree}.

Still it is interesting whether spectral properties correspond for other models.
Since \BA model has triangle-like spectral distribution
\cite{farkas2001spectra}, it should not be compared to \ER model. As the degree
distribution has strong influence on spectral distribution
\cite{nadakuditi_spectra_2013,dorogovtsev2003spectra} the difference may come
from the fact, that \BA model follows power-law degree distribution. Hence \CL
model, which can access arbitrary degree distribution, is better suited model
for spectral comparison.

In this paper we compare \BA model and properly parametrized \CL model in
context of spectral properties. Taking into account the applications in physics
and computer science \cite{anand2011shannon,glos2018vertices,
	chakraborty2016spatial,faccin2013degree}, we focus on the distribution of the
first, the second and the last eigenvalue, the principal eigenvector,
and the spectral distribution of adjacency matrix.

The paper is organized as follows. In Sec.~\ref{sec:preliminaries}, we provide
basic concepts of random graph theory and we specify the research problem. In
Sec.~\ref{sec:numerical-analysis} we specify the experiment and we perform
numerical analysis of spectral similarity between \BA and \CL models. In
Sec.~\ref{sec:conclusion}, we conclude and discuss our results.

\section{Preliminaries}\label{sec:preliminaries}


\BA model $\mathcal G^{\rm BA}_{n}(m_0)$ is an iterative random graph model.
Starting with complete graph of order $m_0$, new vertices are added and
connected randomly with already existing vertices. Instead of connecting them
uniformly at random, vertex $\bar v $ is new neighbour to a newly added vertex with
probability $p_{\bar v} = d_{\bar v}/\sum_vd_v$, where $d_v$ is the degree of
vertex $v$. This way complex graphs are most likely to be generated. Expectedly,
the degree distribution is a power-law distribution of the
form\cite{bollobasdegree}
\begin{equation}
d(k) = \frac{2m_0(m_0+1)}{k(k+1)(k+2)},
\end{equation}
which implies $d(k)\propto k^{-3}$. Contrary to \ER model, the spectrum is not a
semicircle distribution, but mostly it is triangle like
\cite{farkas2001spectra}. Currently, the formula for spectral distribution is
not known.

Contrary to \BA model, in \CL model $\mathcal{G}^{\rm CL}_n(\w)$ the edges appear
independently. The model is a generalization of \ER model with not
necessarily homogeneous degree distribution. Suppose we have real valued vector
$\w=(w_1,\dots,w_n)$ such that $0\leq w_i\leq n-1$. Then, the edge between
vertex $v_i$ and $v_j$ is added with probability
$w_iw_j/\sum_{k=1}^n w_k$ independently. Note that $w_i$ refers to the expected
degree of the vertex $v_i$. The \ER model is recovered for constant $\w$ vector.
Note that \CL is quite a general model, since there are no further restrictions
on $\w$ vector.

The \BA and \CL models are necessarily different. As an example one can consider
the graphs size: for \BA model the size is fixed for fixed $n$ and $m_0$, while
for \CL model the size is a random variable. However despite this fact, it is
possible they have similar spectral properties. More precisely the question is
whether for given value $m_0$ and naturally chosen parameter $\w$, some
spectral properties of \BA and \CL are similar.

In this paper we analyze following spectral properties of adjacency matrix:
\begin{enumerate}
	\item spectral distribution,
	\item distribution of largest eigenvalue,
	\item distribution of second largest eigenvalue,
	\item distribution of smallest eigenvalue,
	\item eigenvector corresponding to the largest eigenvalue (principal eigenvector).
\end{enumerate}
All of them have application in computer science and physics
\cite{anand2011shannon,glos2018vertices,
	chakraborty2016spatial,faccin2013degree}.

Our numerical analysis can be divided into two steps. In the first we determined
$\w$ for each $m_0=1,\dots,6$, which are typical for \BA analysis. In the second
step, for various graph order we determined the values of extreme eigenvalues,
principal eigenvector and spectral distribution. We compared both models using
various similarity measures. For eigenvalues and spectral distribution we chose
the $p$-values of Kolmogorov-Smirnov test. For principal vector we chose the
measures based on uniform norm distance and Euclidean distance.

\section{Numerical analysis} \label{sec:numerical-analysis} 

In this section we analyze numerically the similarity of spectral properties
listed in Sec.~\ref{sec:preliminaries} for \BA and \CL. In
Sec.~\ref{sec:experiment} we describe how the parameter $\w$ for \CL is derived
for given parameter $m_0$ of \BA. In Sec.~\ref{sec:spectrum-bulk} we analyze
similarity of spectral distributions. In Sec.~\ref{sec:extreme-eigenvalues} we
analyze similarity of first, second and last eigenvalue. In
Sec.~\ref{sec:eigenvector} we analyze similarity in principal eigenvectors.

\subsection{The experiment} \label{sec:experiment} Since \CL model has richer
parameter values space, $\w$ parameter is derived on graphs generated by \BA
model. Suppose we fix \BA parameter $m_0$ and the order of graph $n$. Let
$\varepsilon$ by an acceptable deviation for expected degree vector. We start by
generating 300 graphs according to $\mathcal{G}^{\rm BA}_n(m_0)$ and derive
degree vectors $\w^{\rm BA}_{j}$ for $j=1,\dots,300$. We calculate mean degree
vector
\begin{equation}
\overline \w ^{\rm BA}_{1} = \frac{1}{300}\sum_{j=1}^{300}\w^{\rm BA}_{j}.
\end{equation}
Then we generate another 300 graphs and calculate
\begin{equation}
\overline \w ^{\rm BA}_{2} = \frac{1}{600}\sum_{j=1}^{600}\w^{\rm BA}_{j},
\end{equation}
and in general
\begin{equation}
\overline \w ^{\rm BA}_{t} = \frac{1}{300t}\sum_{j=1}^{300t}\w^{\rm BA}_{j}.
\end{equation}
We repeat the process until
\begin{equation}
\| \overline \w ^{\rm BA}_{t} -\overline \w ^{\rm BA}_{t+1} \|_{\infty} \leq \varepsilon.
\end{equation}
We choose $\overline \w ^{\rm BA}_{n,m_0} \coloneqq\overline \w ^{\rm BA}_{t+1}$
to be the parameter for \CL model for given $n$ and $m_0$. Hence we compare
models $\mathcal G_n^{\rm BA}(m_0)$ to $\mathcal G_n^{\rm CL}(\overline \w ^{\rm
	BA}_{n,m_0})$. In our experiment we considered $m_0 =1,\dots,6$ and
$\varepsilon=0.05$. The values of graph orders depends on analyzed spectral
property.

The numerical analysis was performed using Julia language. The graphs were
generated thanks to \texttt{LightGraphs.jl} module \cite{Bromberger17}. The
module implemented the algorithm for \CL generation presented in
\cite{miller2011efficient}. The algorithm scales well for expected degrees much
smaller than the graph order which is true for \BA model \cite{flaxman2005high}

\subsection{Spectral distribution}\label{sec:spectrum-bulk}

\begin{figure}[t]\centering
	\begin{tabular}{@{}l@{\!\!\!}l@{\!\!\!\!}l@{}}
		\includegraphics{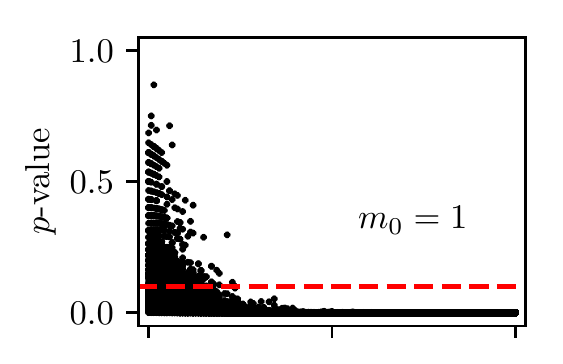}&
		\includegraphics{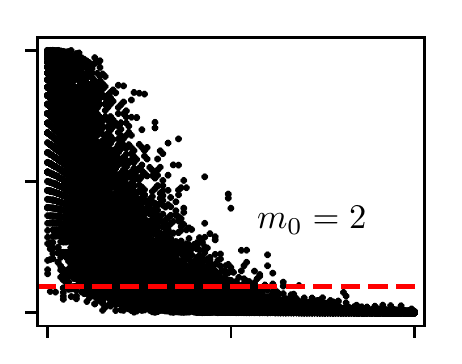}&
		\includegraphics{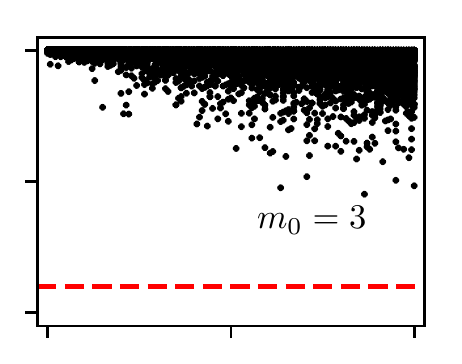}\\[0em]
		\includegraphics{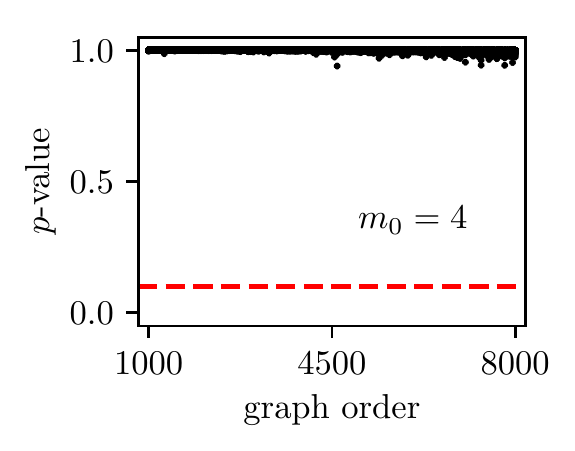}&
		\includegraphics{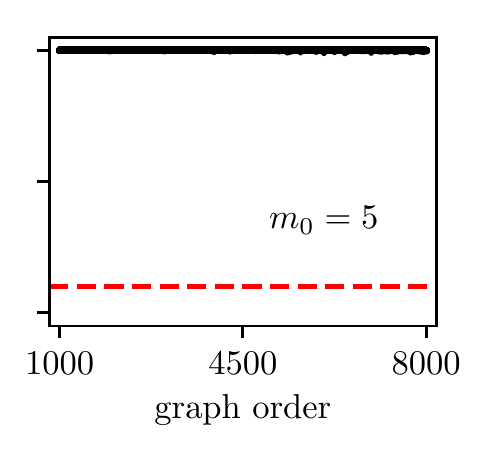}&
		\includegraphics{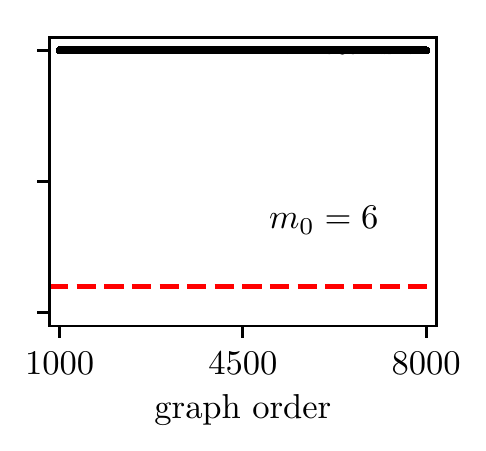}
	\end{tabular}
	\caption{$p$-value of Kolmogorov-Smirnov distribution equality test. Red dashed
	lines refer to critical value $\alpha=0.1$.}\label{fig:spectrum-bulk}
\end{figure}

The experiment goes as follow. For each order $n=1000,1050,\dots,8000$ and
$m_0=1,\dots,6$ we have calculated the spectrum of 200 graphs for each model. We
have paired consecutive spectra and perform a Kolmogorov-Smirnov test for each
pair. By this for each order and $m_0$ value we obtained 200 $p$-values. Our
results are presented in Fig.~\ref{fig:spectrum-bulk}.

We can see, that for $m_0=1,2$ the spectral distributions differ, as for large
graph order the $p$-values converge to 0. While for $m_0=3$ the $p$-value are
close to 1 for small graph orders, we can see that for large $n$ the $p$-values
are more robust. For $m_0>3$ almost all $p$-values are very close to 1.

The $m_0=3$ is a threshold value for spectral distribution similarity between
the models. Below that value, the spectra are almost surely different. We cannot
specify any reliable statement for $m_0=3$, however it seems for sufficiently
`small' graphs the distributions are similar. For larger values of parameters
distributions are similar.

The possible reason of such behavior is that for small values of $m_0$ we obtain
in the case of Chung-Lu model large number of disconnected small graphs. The
greater the $m_0$ is, the smaller is the number of such components. Above the
threshold value $m_0$ disconnected components have marginal impact on the
spectrum hence the models spectra are similar. Because of that we analyzed
principal vectors and extreme eigenvalues only for $m_0=4,5,6$.

\subsection{Extreme eigenvalues}\label{sec:extreme-eigenvalues}

\begin{figure}[t]
	\begin{tabular}{@{}l@{\!\!\!}l@{\!\!\!\!}l@{}}
		\includegraphics{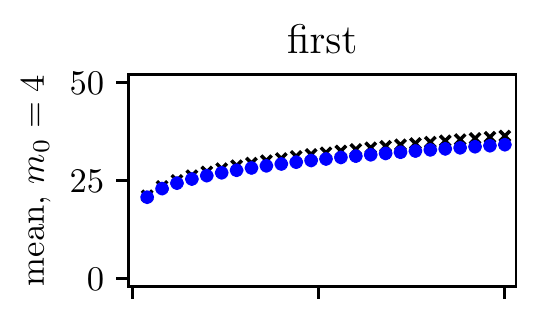}&
		\includegraphics{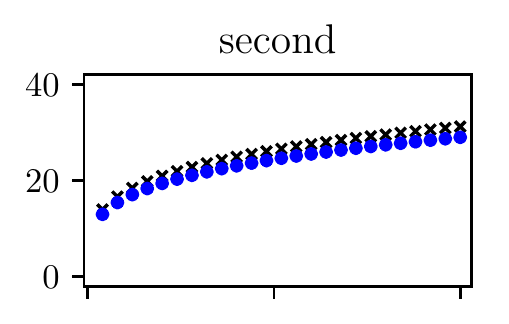}&
		\includegraphics{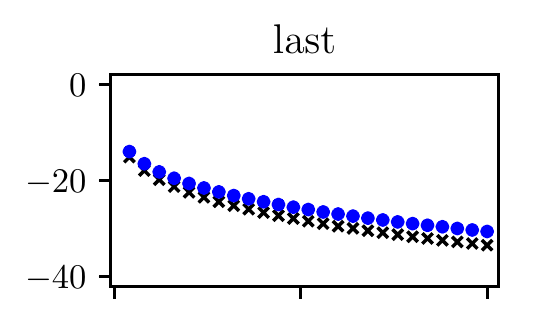} \\[0em]
		\includegraphics{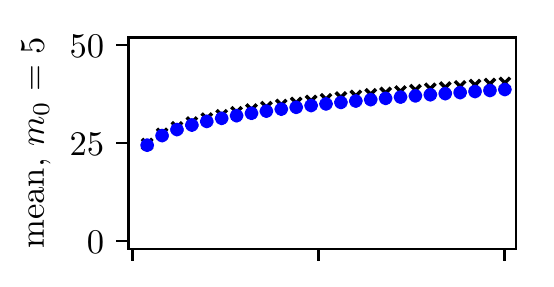}&
		\includegraphics{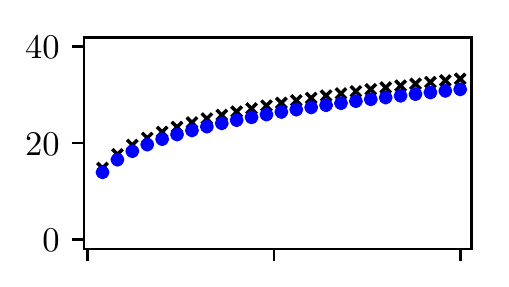}&
		\includegraphics{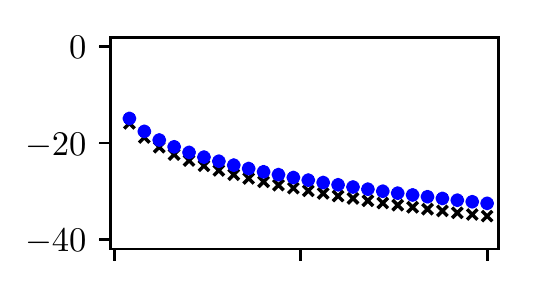} \\[0em]
		\includegraphics{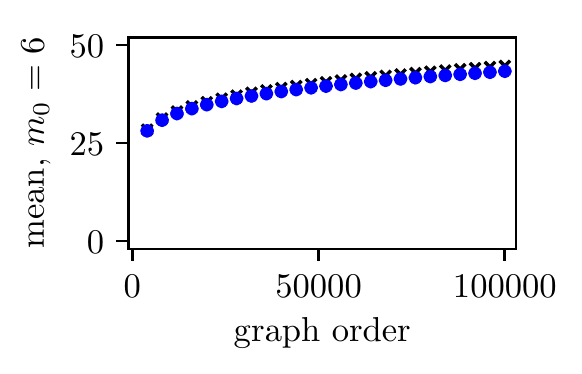}&
		\includegraphics{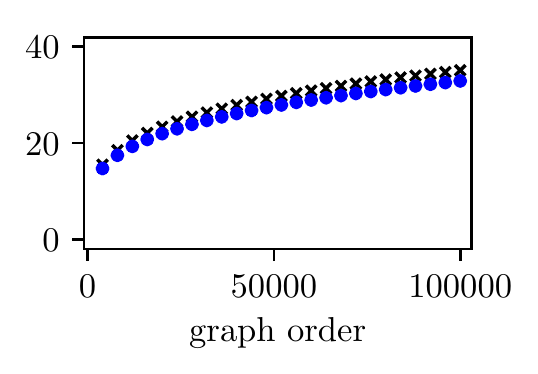}&
		\includegraphics{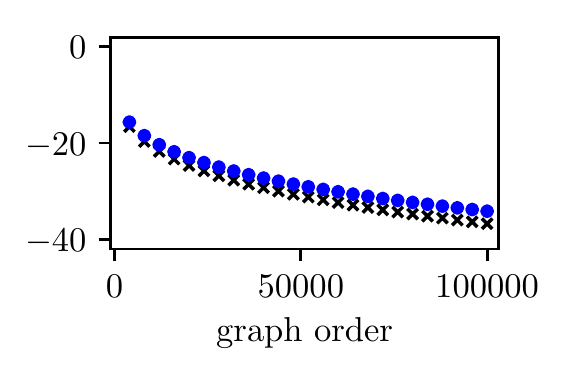}
		
	\end{tabular}
\caption{Mean value of first, second and last eigenvalues. The black crosses
	correspond to \BA model, while blue dots correspond to \CL model. }\label{fig:eigenvalues-mean}
\end{figure}

\begin{figure}[t]
	\begin{tabular}{@{}l@{\!\!\!}l@{\!\!\!\!}l@{}}
		\includegraphics{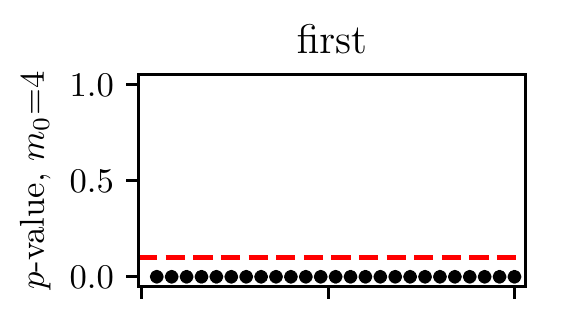}&
		\includegraphics{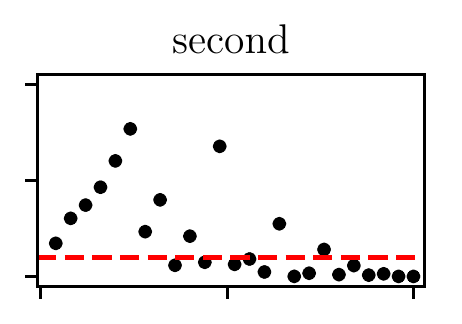}&
		\includegraphics{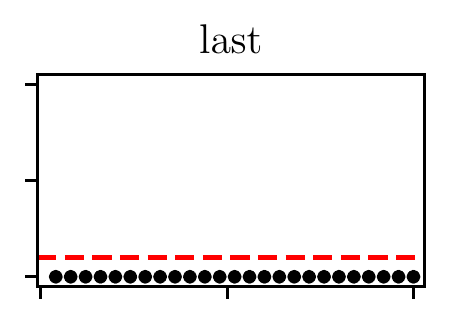}\\[0em]
		\includegraphics{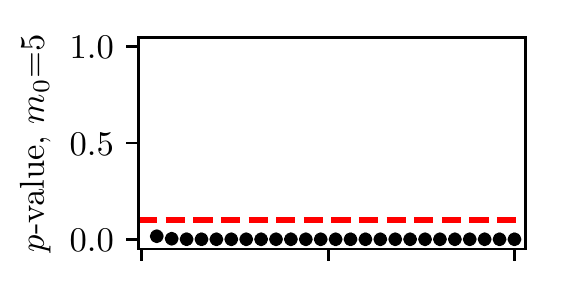}&
		\includegraphics{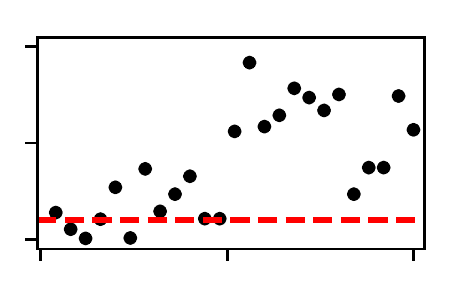}&
		\includegraphics{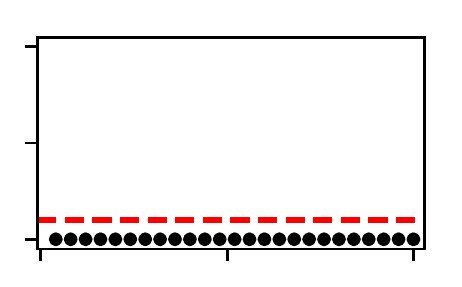}\\[0em]
		\includegraphics{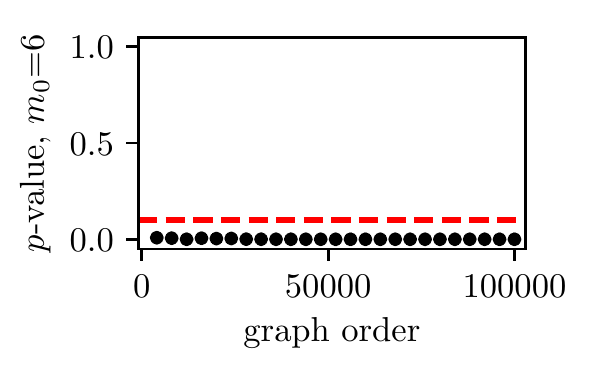}&
		\includegraphics{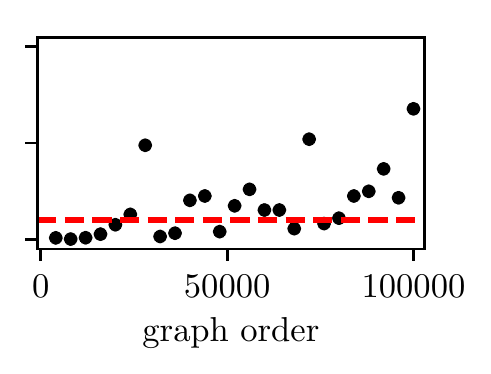}&
		\includegraphics{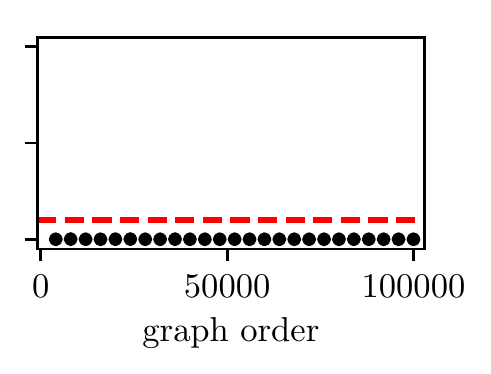}
	\end{tabular}
	\caption{$p$-values of standardized eigenvalues data. Red dashed lines refer to
	critical value $\alpha=0.1$.}\label{fig:eigenvalues-pvalue}
\end{figure}

The experiment goes as follow. For each graph order
$n=4000,8000,\dots,100\,000$ and $m_0=4,5,6$ we have calculated first, second
and last eigenvalues of $20,\!000$ graphs. Then we compare empirical
distributions of \BA and \CL model by calculating mean value, and $p$-value of
Kolmogorov-Smirnov for standardized data. By this for each graph order and $m_0$
value we obtained single mean value and $p$-value. Our results are presented in
Fig.~\ref{fig:eigenvalues-mean} and Fig.~\ref{fig:eigenvalues-pvalue}.

As one can see on the Fig.~\ref{fig:eigenvalues-mean} the mean values slightly
differ. While the gap between mean values increases with the graph order growth,
the values are still quite similar. However $p$-values presented in
Fig.~\ref{fig:eigenvalues-pvalue} show that independently of $m_0$ distributions
of first and last eigenvalues differ. For $m_0\leq 4$ this is the case for
second eigenvalue.

Surprisingly, we cannot provide any reliable statement for second eigenvalue for
$m_0=5,6$, as $p$-values appears on both side of $\alpha=0.1$. While most of the $p$-values are above
critical value $\alpha=0.1$, the result are not as evident as in
Sec.~\ref{sec:spectrum-bulk}. Hence we cannot decline the possibility, that
those distributions are similar.
\subsection{Principal eigenvector}\label{sec:eigenvector}

While $p$-value  was a natural choice for eigenvalues and spectral distribution analysis,
this is not the case for eigenvector. For graph of order $n$, $n$-dimensional
variable would be considered. Furthermore, since the length of the vector goes
to infinity, norms cease to be equivalent. Hence, we chose the following measures of
similarity between vectors $\bar \lambda,\bar \kappa$:
\begin{gather}
	\|\bar \lambda -\bar \kappa \|_\infty = \max_{i=1,\dots,n} |\bar{\lambda}_i-\bar{\kappa}_i|,\\
	\|\bar \lambda -\bar \kappa \|_2 = \sqrt{\frac{1}{2}\sum_{i=1}^{n} (\bar{\lambda}_i-\bar{\kappa}_i)^2}.
\end{gather}
The constant $\frac{1}{\sqrt{2}}$ was chosen to bound the possible values to
$[0,1]$. Note that $\|x\|_2\to 0$ implies $\|x\|_\infty \to
0$, hence the first measure is the stronger one.

The choice of eigenvector is in general non-unique, as for
eigenvector $\bar \lambda$ and arbitrary $\varphi\in\R$, vector  $\ee^{\ii\varphi}\bar
\lambda_1$ is proper eigenvector corresponding to the same eigenvalue. However
based on Perron-Frobenius theorem, we choose the unique eigenvector with all
elements nonnegative. Note that in the case of similarity both measures should
converge to 0. 

\begin{figure}[t]\centering
	\begin{tabular}{c@{\hspace{-0.1cm}}c}
		\includegraphics{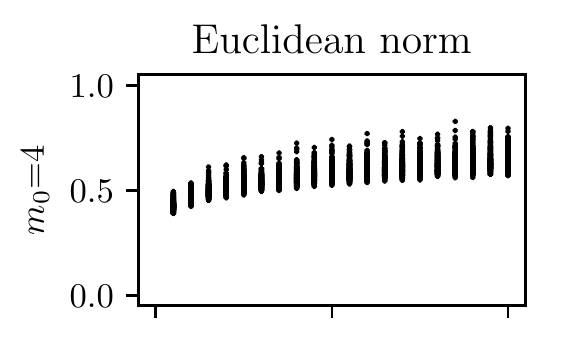}&
		\includegraphics{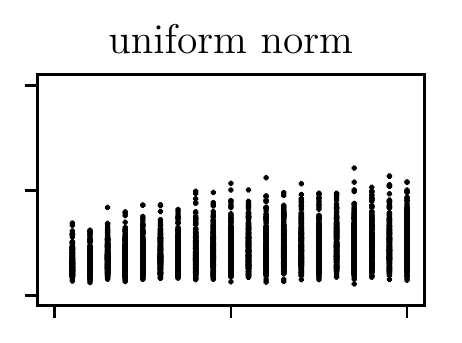}\\[0em]
		\includegraphics{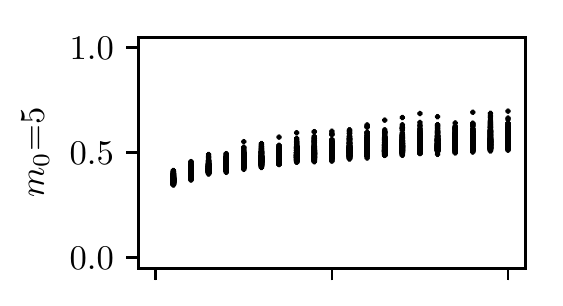}&
		\includegraphics{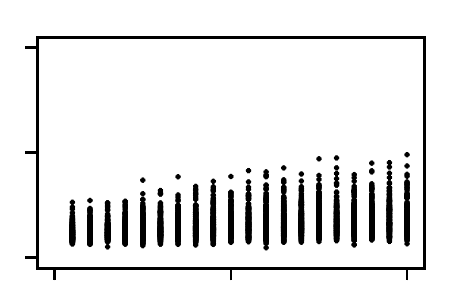}\\[0em]
		\includegraphics{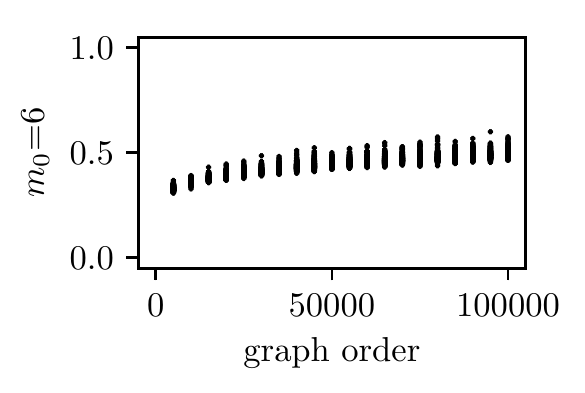}&
		\includegraphics{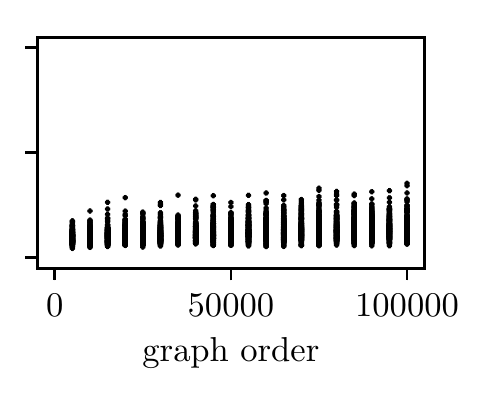}
	\end{tabular}
	\caption{Various similarity measures of principals vector for \BA and \CL
	model. In the case of similarity, both measures should converge to
	0.}\label{fig:eigenvector}
\end{figure}

Based on the numerical results presented in Fig.~\ref{fig:eigenvector} this is
not the case. For all $m_0$ values Euclidean distance $\|\cdot\|_2$ is detached from 0.
In the case of infinity norm distance, we observe more robust behavior. There
exist values close to 0, however still for many samples the values are further
from 0. Hence in our opinion it is unlikely for them to be equal as well.

\section{Application of the results} To demonstrate the possible applications of
our results we consider two examples. In the first one we show that spectral
distribution similarity can be used for derivation of the spectral distribution
formula for \BA model. In the second case we show that the similarities and
differences of eigenvalues and eigenvectors have impact on the quantum spatial
search evolution.

\subsection{Spectral distribution}

\begin{figure}[t]\centering
	\begin{tabular}{@{}l@{\!\!\!}l@{\!\!\!\!}l@{}}
		\includegraphics{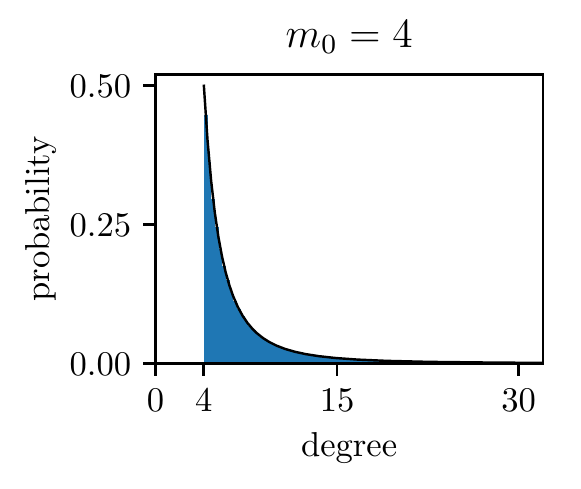}&
		\includegraphics{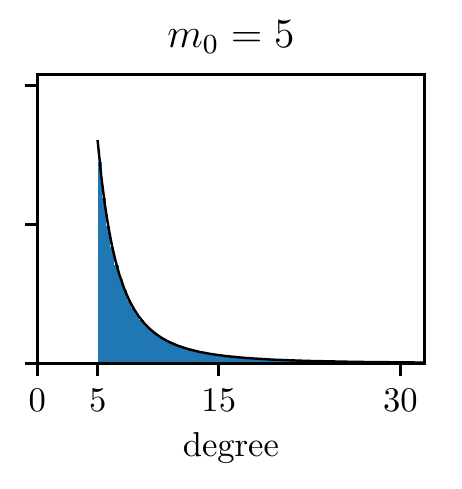}&
		\includegraphics{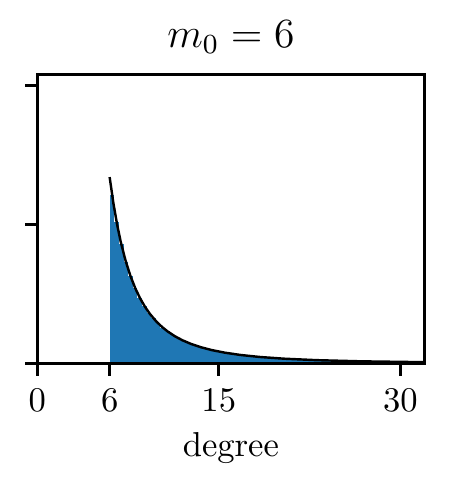}
	\end{tabular}
	\caption{Expected degree distribution of \BA model. The black line refers the
	analytical formula $2m_0^2/d^3$. The data were constructed as described in
	Sec.~\ref{sec:experiment} for graph order $100,000$.}\label{fig:degree}
\end{figure}

While spectral distributions are known for some models
\cite{nadakuditi_spectra_2013,dorogovtsev2003spectra}, it is not the case of \BA
model. Thanks to our results obtained in Sec.~\ref{sec:spectrum-bulk}, we can
reduce the problem of deriving it to derivation of much simpler Chung-Lu model.
We propose following approach: First, one should shown that the dependence
between edges appearance has no influence on the \BA spectral distribution,
which is numerically confirmed in the paper. Second, using methods suitable for
edge-independent random graphs we may determine the analytical formula for
spectral distribution. The first method requires degree distribution, which is
well known for the \BA model \cite{bollobasdegree}.  The methods proposed by
Dorogovstev et al. \cite{dorogovtsev2003spectra} allow to use this
characteristics for spectral distribution derivation.

The second method proposed by Nadakuditi et al. \cite{nadakuditi_spectra_2013}
requires the distribution of the expected degree. While it is not known for \BA
model, we propose $p:[m_0,\infty)$ such that $p:d\mapsto\frac{2m_0^2}{d^3}$
distribution. The formula fits the histogram very well, furthermore it is well
motivated because of several analytical requirements. First, it is decreasing
function which scales as $p \propto d^{-3}$ \cite{albert2002statistical}.
Second, $\mathbb E p(d) = 2m_0$, which coincide with the mean degree for the
model.

\subsection{Quantum spatial search} \label{sec:qss-application}

\begin{figure}[t]\centering
	\begin{tabular}{@{}l@{\!\!\!}l@{\!\!\!\!}l@{}}
		\includegraphics{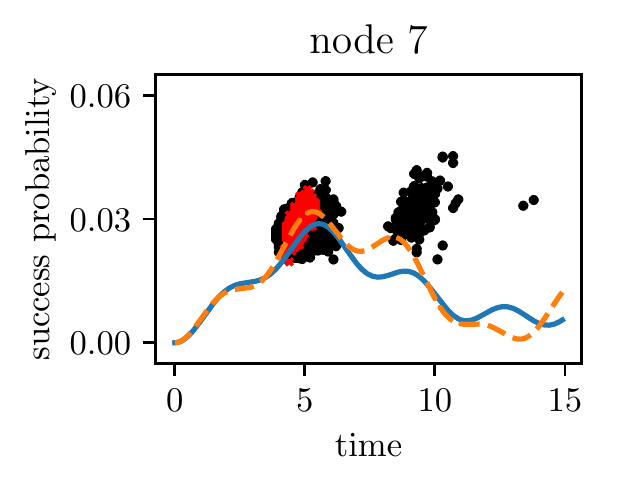}&
		\includegraphics{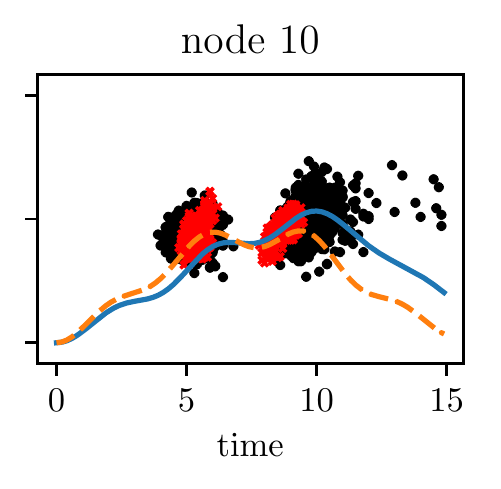}&
		\includegraphics{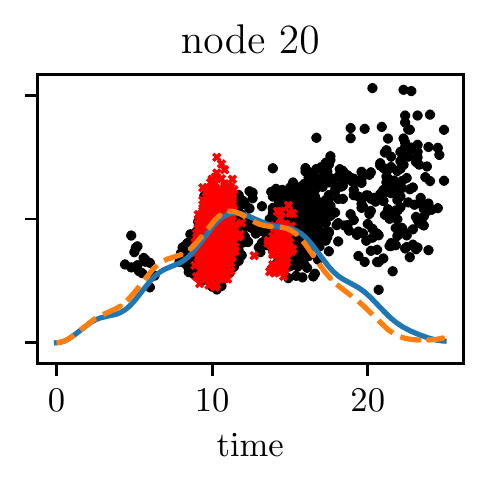}
	\end{tabular}
	\caption{The optimal measurement times and probability success for \BA and \CL
	models. Dashed orange represent a success probability for single \CL graph.
	Blue line represent a success probability change for single \BA graph. Red
	crosses (black dots) represent an optimal measurement time for \CL graphs (\BA
	graphs). The procedure of deriving optimal measurement time are described in
	Sec.~\ref{sec:qss-application}  }\label{fig:qss}
\end{figure}

\begin{figure}[t]\centering
	\begin{tabular}{@{}l@{\!\!\!}l@{\!\!\!\!}@{}}
		\includegraphics{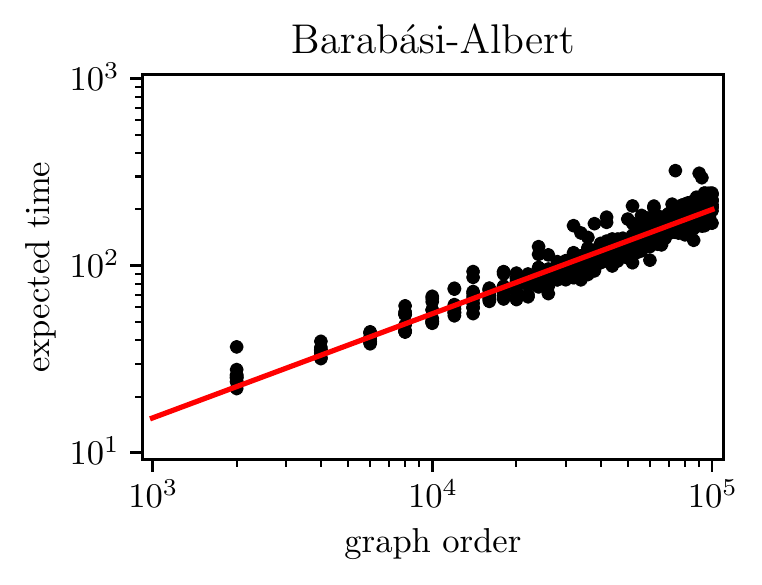} &
		\includegraphics{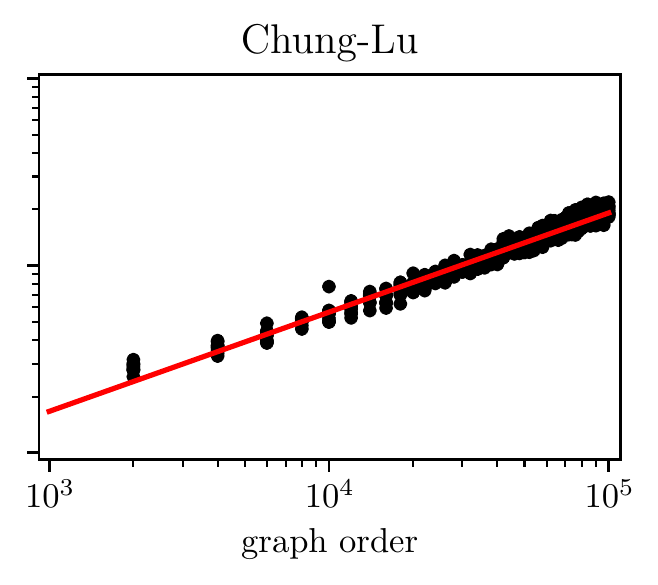}
	\end{tabular}
	\caption{The minimal expected time $(t+0.2\log(n))/p(t)$ for \BA and \CL model. The red line is a regression line of logarithm of expected time vs logarithm of graph order. It describes the time complexity of the quantum search, as tangent $\alpha$ implies the complexity is $\Theta(n^\alpha)$. For \BA model its tangent equals $\alpha\approx0.56$, while for \CL it is $\alpha\approx 0.53$}\label{fig:qss-scaling}
\end{figure}

The second example is motivated by the fact, that extreme eigenvalues and
principal eigenvector play important role in quantum spatial search efficiency
\cite{chakraborty2016spatial,glos2018vertices}. Let
\begin{equation}
H = -t(\gamma A + \proj{\omega})
\end{equation}
be a Hamiltonian describing continuous-time quantum spatial search
\cite{childs2004spatial}. Here $A$ denotes adjacency matrix of a graph and
$\proj{\omega}$ is an oracle of marked vertex, and $t$ denotes evolution time.
Jumping rate $\gamma$ is chosen to be $\frac{1}{\lambda_1(A)}$. After time $t$ a
measurement in computational basis is made in order to derive $\omega$. It was
shown, that large gap between largest eigenvalue and the maximum of absolute
values of second and last eigenvalue implies optimal $\sqrt n$ time of quantum
evolution \cite{chakraborty2016spatial}. While the result applies only to
regular graphs, the results from \cite{glos2018vertices} suggests this may
depend mostly on the overlap between the principal eigenvector and the subspace
corresponding to marked element.

We have designed experiment as follows. We made an evolution of quantum spatial
search on 1800 graphs from $\mathcal G_{100\,000}^{BA}(6)$ and 1800 graphs from 
$\mathcal G_{100\,000}^{CL}(\bar {\mathbf{w}}^{BA}_{100\,000,6})$. The evolution
was made for marked nodes 7, 10 and 20. For marked nodes 7 and 10 we determined
the success probability for times $0,0.1,\dots,15$, while for marked node 20 we
have chosen times $0,0.1\dots, 25$. Then we chose the time, for which the
success probability $p_{\rm opt}$ is the first greatest probability, for which
all of the probabilities $p_t$ for larger evolution time satisfies $(p_t-p_{\rm
	opt})/p_t <0.2$. This prevents from choosing success probability too far from $t=0$,
which would result in not necessarily optimal quantum spatial search evolution.
For simulation purposes we used the Julia language together with \texttt{QuantumWalk.jl} package \cite{glos2018quantumwalkjl}.

The results are presented in Fig.~\ref{fig:qss}. We can observe that there are
regimes, in which the success probabilities are similar for both models for each
marked nodes. However, the greater the index of node is, the more robust results
are observed. We claim that the observed behavior can be explained by the
difference in the principal eigenvector convergence presented in
Sec.~\ref{sec:eigenvector}.

The results is stable with graph order change. For each graph order $n=2000, 4000,\ldots, 100,000$ we
have sampled 10 graphs for each model and calculated the optimal expected time for
quantum search with marked $(m_0+1)$th node. We have calculated the success probability
$p(t)$ for times $t=0,0.01\pi\sqrt{n},\ldots,\pi\sqrt{n}$ and we have chosen
time $t_{\rm opt}$ minimizing the expected time
$(t_{\rm opt} +0.1\log(n))/p(t)$. The approach was used already for analyzing
the quantum attack \cite{glos2018impact}. The results are presented in
Fig.~\ref{fig:qss-scaling}. We can see that the optimal expected time is similar
for both models. The same effect can be observed for scaling, as the calculated
complexity is $\Theta(n^{0.56})$ for \BA model and $\Theta(n^{0.53})$ for \CL
model. This confirms that the similarity of the quantum spatial search behavior
scales with graphs order.

\section{Conclusion and discussion} \label{sec:conclusion}

In the paper we analyzed the spectral similarity between \BA model and \CL
model. Our analysis includes the similarity of distribution of spectral
distribution, extreme eigenvalues and principal eigenvector.

Calculations show that $m_0=3$ is threshold value for the spectral
distributions similarity. The models may have similar spectral distributions for
$m_0>3$ and different distributions for $m_0<3$. Contrary the extreme
eigenvalues seems not to be similar. In particular for all considered values of
$m_0$, first and last eigenvalues of those models differ. Still we cannot
provide any reliable hypothesis for $m_0=5,6$.

For principal vectors we have proposed Euclidean and uniform norm distance
similarity measures. In the first case for analyzed values of $m_0$ the
principal vector differs, as the distance was detached from 0. We cannot provide
such strong statement for the infinity norm measure, due to robustness of
obtained values. We have made similar analysis for extreme eigenvalues, here the
eigenvalues are roughly the same. Still the distributions differ at least for
largest and smallest eigenvalues.

We have proposed applications of our results in graph spectral theory and
quantum spatial search. In the first case we proposed a method for deriving
analytical formula for spectral distribution of quantum spatial search. In the
second we have shown, that robustness of convergence of principal eigenvector
has impact on the quantum spatial search optimal measurement time. Still, for
 vertices with large degree the results are similar for both models. As a future research one can consider extending the results into other quantum walks problems like quantum transport \cite{mulken2011continuous}.

\paragraph{Acknowledgements} The author would like to thank Łukasz Pawela,
Zbigniew Puchała and Jarosław Adam Miszczak for discussion. Furthermore the
author would like to thank Izabela Miszczak and Jarosław Adam Miszczak for
reviewing the article. This research was supported in part by PL-Grid
Infrastructure. The work was carried out within the statutory research project of the Institute of Informatics, BKM.

\bibliographystyle{ieeetr}
\bibliography{random-graph-spectrum-comparison}

\end{document}